\documentclass[preprint,timesc]{elsarticle}

\usepackage{amssymb,amsmath,amsthm}
\usepackage{bm}
\usepackage{booktabs}
\usepackage{algorithm}
\usepackage{algpseudocode}
\usepackage{threeparttable}
\usepackage{geometry}
\usepackage[T1]{fontenc}  
\usepackage{newpxtext, newpxmath}  

\newtheorem{theorem}{Theorem}
\newtheorem{lemma}{Lemma}

\newtheorem{corollary}{Corollary}
\theoremstyle{definition}

\newtheorem{example}{Example}
\theoremstyle{remark}
\newtheorem{remark}{Remark}
\newtheorem{assumption}{Assumption}

\newcommand{\R}{\mathbb{R}}
\newcommand{\Z}{\mathbb{Z}}
\newcommand{\X}{\mathcal{X}}
\newcommand{\vt}[1]{\boldsymbol{#1}}
\newcommand{\mini}{\text{minimize}}
\newcommand{\st}{\mbox{subject to}}
\DeclareMathOperator{\median}{median}
\DeclareMathOperator{\OPT}{OPT}

\allowdisplaybreaks

\journal{Elsevier}

\begin{document}
\begin{frontmatter}

\title{%
    Inverse optimization for divisor methods of apportionment
}

\author[jgmi]{Daiki Kawane}
\ead{kawane.daiki.400@s.kyushu-u.ac.jp}

\author[econ_grad]{Hokuto Nagano}
\ead{nagano.hokuto.397@s.kyushu-u.ac.jp}

\author[econ]{Yoichi Izunaga\corref{cor}}
\ead{izunaga@econ.kyushu-u.ac.jp}

\address[jgmi]{%
    Joint Graduate School of Mathematics for Innovation,
    Kyushu University,
    Fukuoka,
    Japan
}

\address[econ_grad]{%
    Graduate School of Economics,
    Kyushu University,
    Fukuoka,
    Japan
}

\address[econ]{%
    Faculty of Economics,
    Kyushu University,
    Fukuoka,
    Japan
}

\cortext[cor]{Corresponding author.}

\begin{abstract}
The apportionment problem is to allocate a fixed number of seats among regions in proportion to their populations.
In this paper, we propose an inverse optimization approach to evaluate how fairly an observed apportionment treats each region.
It is well known that every divisor method can be formulated as a separable discrete convex minimization problem.
We introduce a weighted divisor method, in which a nonnegative weight is assigned to each region, as the forward problem.
Then we formulate an inverse optimization problem that finds the weight vector closest to the all-one vector in the $\ell_1$-norm such that the observed apportionment is optimal for the forward problem.
We show that the optimal weight vector is always unique and has a closed-form expression as the median of three values that depend on a single parameter.
Based on this result, the inverse optimization problem is reduced to a one-dimensional monotone equation, and it can be solved in $O(n \log n)$ time, where $n$ is the number of regions.
We also apply the proposed approach to the apportionment of the single-member districts of the House of Representatives of Japan.
\end{abstract}

\end{frontmatter}

\section{Introduction}
\label{sec:intro}

The apportionment problem is to allocate a fixed number of seats among regions, such as states or prefectures, in proportion to their populations~\cite{BalinskiYoung2001,Pukelsheim2017}.
Since the number of seats must be an integer, the exact proportional allocation cannot be realized in general, and some rounding rule is necessary.
Divisor methods are one of the most important classes of apportionment methods.
For example, the Hill method (also called the Huntington--Hill method) has been used for the House of Representatives of the United States, and the Adams method has been adopted for the House of Representatives of Japan.

When an apportionment is actually decided, it is important to explain quantitatively how fairly the apportionment treats each region.
In the House of Representatives of Japan, the seats of the single-member districts had been apportioned among the 47 prefectures by the so-called \emph{one-seat-per-prefecture rule}, in which each prefecture first receives one seat and the remaining seats are allocated in proportion to the population.
Since this rule favors prefectures with small populations, the disparity in the value of a vote has been a serious problem, and the apportionment has been revised several times.
However, it is not easy to evaluate such revised apportionments, since they are not always obtained by a specific apportionment method.

In this paper, we evaluate an observed apportionment by using \emph{inverse optimization}.
In a usual (forward) optimization problem, the parameters are given and we find an optimal solution.
In an inverse optimization problem, a solution is observed, and we find parameters that make the observed solution optimal, with the minimum change from given nominal parameters~\cite{AhujaOrlin2001,Heuberger2004,ChanMahmoodZhu2025}.
Our approach is based on the well-known fact that every divisor method can be formulated as a separable discrete convex minimization problem under a single cardinality constraint~\cite{IbarakiKatoh1988}.
We introduce a \emph{weighted divisor method} as the forward problem, where a nonnegative weight $w_i$ is assigned to the objective function of each region $i$.
The weight $w_i$ represents the relative priority of region~$i$: a region with $w_i<1$ obtains seats more easily, and a region with $w_i>1$ obtains seats less easily.
Then we consider the inverse optimization problem that finds the weight vector closest to the all-one vector $\vt{1}$ in the $\ell_1$-norm such that the observed apportionment is optimal for the weighted divisor method.
The deviation of the optimal weights from $\vt{1}$ can be used as a measure of fairness of the observed apportionment with respect to the chosen divisor method.

The main contributions of this paper are as follows.
\begin{itemize}
  \item For a general divisor method, we show that the optimal weight vector of the inverse optimization problem is always unique, and that each component is expressed as the median of three values depending on a one-dimensional parameter $\lambda$ (Theorem~\ref{thm:median}).
  \item By this median representation, the inverse optimization problem is reduced to a one-dimensional monotone equation $S(\lambda)=n$.
  This is analogous to the forward problem, in which the divisor $r$ is adjusted so that the total number of seats is equal to $K$.
  Based on this reduction, we give an $O(n\log n)$ time algorithm (Algorithm~\ref{alg:IP}).
  \item We show an exchange property of the optimal weights (Theorem~\ref{thm:exchange}): if one seat is moved from a region with $w^*_i<1$ to a region with $w^*_j>1$, then the observed apportionment becomes closer to the one obtained by the chosen divisor method. Thus, the optimal weights indicate not only the degree of unfairness but also a direction of correction.
\end{itemize}

The rest of this paper is organized as follows.
In Section~\ref{sec:prelim}, we review divisor methods and their formulation as a separable discrete convex minimization problem.
In Section~\ref{sec:inverse}, we introduce the weighted divisor method and formulate the inverse optimization problem.
In Section~\ref{sec:theory}, we show the median representation, the uniqueness, and an exchange property of the optimal weights, and give an algorithm based on binary search.
Section~\ref{sec:experiment} applies the proposed approach to the House of Representatives of Japan, and Section~\ref{sec:conclusion} concludes the paper.

\section{Preliminaries}
\label{sec:prelim}

\subsection{Divisor methods}

Throughout this paper, for an integer $k$, we denote by $\Z_{\ge k}$ the set of integers greater than or equal to $k$, and by $\R_{\ge 0}$ the set of nonnegative real numbers.
Let $N=\{1,2,\dots,n\}$ be the set of regions, and let $p_i>0$ be the population of region $i\in N$.
Let $K$ be the total number of seats; we assume $K\ge n$.
An \emph{apportionment} is a vector $\vt{x}=(x_1,\dots,x_n)^\top$ in $\X = \{ \vt{x}\in\Z_{\ge 0}^n \mid \sum_{i\in N} x_i = K \}$,
where $x_i$ is the number of seats of region $i$.
If we regard proportionality to population as the criterion of fairness, the ideal number of seats of region~$i$ is the \emph{quota} $q_i = K p_i / \sum_{j\in N} p_j$.
However, $q_i$ is not an integer in general, and thus some rounding procedure is required.

A \emph{divisor criterion} is a function $d\colon \Z_{\ge 0}\to\R$ satisfying
\begin{align}
  k \le d(k) \le k+1 \quad\text{and}\quad d(k) < d(k+1) \qquad (k\in\Z_{\ge 0}).
  \label{eq:divisor-criterion}
\end{align}
For notational convenience, we set $d(-1):=0$.
For a real number $q\ge 0$, the \emph{$d$-rounding} of $q$ is the set
\begin{align*}
  [q]_d := \{\, x\in\Z_{\ge 0} \mid d(x-1)\le q \le d(x) \,\}.
\end{align*}
The \emph{divisor method} based on $d$ determines an apportionment as follows.
\begin{enumerate}
  \item Choose an initial divisor $r>0$.
  \item For each $i \in N$, compute the quotient $p_i/r$ and round it as $x_i(r) \in [p_i/r]_d$.
  \item If $\sum_{i\in N} x_i(r) = K$, then stop. Otherwise, adjust $r$ and go to Step~2.
\end{enumerate}
In other words, $\vt{x}\in\Z_{\ge 0}^n$ is an apportionment obtained by the divisor method based on $d$ if and only if there exists $r>0$ such that $x_i\in[p_i/r]_d$ for all $i\in N$ and $\sum_{i\in N}x_i=K$.
Note that $[q]_d$ consists of two integers $x$ and $x+1$ if $q=d(x)$ for some
$x\in\Z_{\ge0}$; in this case, all apportionments satisfying the above condition are
regarded as being obtained by the divisor method.

Five representative divisor methods are the Adams method ($d(x)=x$, rounding up), the Dean method ($d(x)=x(x+1)/(x+1/2)$, harmonic mean), the Hill method ($d(x)=\sqrt{x(x+1)}$, geometric mean), the Webster method ($d(x)=x+1/2$, standard rounding), and the Jefferson method ($d(x)=x+1$, rounding down).
They satisfy $d_{\mathrm{Adams}}\le d_{\mathrm{Dean}}\le d_{\mathrm{Hill}}\le d_{\mathrm{Webster}}\le d_{\mathrm{Jefferson}}$.
If $d(0)=0$, as in the Adams, Dean, and Hill methods, then no region receives zero seats.
Roughly speaking, a smaller $d$ rounds up quotients more often.
Since a region with a small population has a small quotient $p_i/r$, one additional seat by rounding up has a relatively large effect for such a region.
Therefore, a divisor method with a smaller $d$ tends to favor regions with small populations~\cite{BalinskiYoung2001}.

\subsection{Formulation as separable discrete convex minimization}

For each $i\in N$, define the function $f_i\colon\Z_{\ge 0}\to\R$ by $f_i(x_i) = (1/p_{i})\sum_{k=0}^{x_i-1} d(k)$,
where the empty sum is zero, and consider the following problem:
\begin{equation*}
  \text{(P)} \qquad \mini \ \sum_{i\in N} f_i(x_i) \quad \st \,\, \vt{x}\in\X.
\end{equation*}
We define $d_{i}(x)=d(x)/p_{i}$ for each $i \in N$ and $x \in \Z_{\ge -1}$.
Then $f_i(x_i+1)-f_i(x_i) = d_i(x_i)$ holds for $x_i\ge 0$.
Since $d$ is increasing, each $f_i$ is a discrete convex function, that is, $f_i(x+1)-f_i(x)$ is nondecreasing in $x$.
Thus, (P) is a separable discrete convex resource allocation problem.

We recall some notions from discrete convex analysis~\cite{Murota2003}.
For $i\in N$, let $\vt{e}_i$ denote the $i$th unit vector.
A set $S\subseteq\Z^n$ is called \emph{M-convex} if, for any $\vt{x},\vt{y}\in S$ and $i\in N$ with $x_i>y_i$, there exists $j\in N$ with $x_j<y_j$ such that $\vt{x}-\vt{e}_i+\vt{e}_j\in S$ and $\vt{y}+\vt{e}_i-\vt{e}_j\in S$.
A function $g\colon\Z^n\to\R\cup\{+\infty\}$ is called \emph{M-convex} if, for any $\vt{x},\vt{y}$ in its effective domain $\{\vt{z}\in\Z^n\mid g(\vt{z})<+\infty\}$ and $i\in N$ with $x_i>y_i$, there exists $j\in N$ with $x_j<y_j$ such that
\begin{equation*}
  g(\vt{x})+g(\vt{y})\ge g(\vt{x}-\vt{e}_i+\vt{e}_j)+g(\vt{y}+\vt{e}_i-\vt{e}_j).
\end{equation*}
The set $\X$ is M-convex, and a separable discrete convex function restricted to $\X$ (i.e., set to $+\infty$ outside $\X$) is an M-convex function.
Hence, the following optimality condition is a special case of the optimality criterion for M-convex functions~\cite{Murota2003}.

\begin{lemma}[{\cite[Theorem~9.8]{Murota2024}}]\label{lem:optimality}
Let $g_i\colon\Z_{\ge 0}\to\R$ $(i\in N)$ be discrete convex functions, and let $\delta_i(x):=g_i(x+1)-g_i(x)$ for $x\ge 0$.
Then $\vt{x}\in\X$ is an optimal solution of $\min\{\sum_{i\in N} g_i(x_i) \mid \vt{x}\in\X\}$ if and only if
\begin{equation*}
  \delta_i(x_i-1) \le \delta_j(x_j) \qquad \text{for all } i,j\in N \text{ with } x_i\ge 1.
\end{equation*}
\end{lemma}

Lemma~\ref{lem:optimality} says that $\vt{x}$ is optimal if and only if no exchange of one seat from a region to another region decreases the objective value.
By using Lemma~\ref{lem:optimality}, we obtain the following theorem, which shows that (P) represents the divisor method.

\begin{theorem}\label{thm:divisor-P}
An apportionment $\vt{x}\in\X$ is obtained by the divisor method based on $d$ if and only if $\vt{x}$ is an optimal solution of (P).
Moreover, $\vt{x}\in\X$ is optimal for (P) if and only if there exists $\lambda>0$ such that
\begin{align}
  d_i(x_i-1) \le \lambda \le d_i(x_i) \qquad (i\in N),
  \label{eq:opt-P}
\end{align}
and in this case $r=1/\lambda$ is a divisor that yields $\vt{x}$.
\end{theorem}

\begin{proof}
For $\lambda=1/r>0$, the condition $x_i\in[p_i/r]_d$ is equivalent to $d(x_i-1)\le p_i\lambda\le d(x_i)$, that is, $d_i(x_i-1)\le\lambda\le d_i(x_i)$.
Hence, it suffices to show that $\vt{x}\in\X$ is optimal for (P) if and only if \eqref{eq:opt-P} holds for some $\lambda>0$.

Suppose that \eqref{eq:opt-P} holds.
Then $d_i(x_i-1)\le\lambda\le d_j(x_j)$ for all $i,j\in N$, and $\vt{x}$ is optimal by Lemma~\ref{lem:optimality}.

Conversely, suppose that $\vt{x}$ is optimal.
By Lemma~\ref{lem:optimality} and $d_i(-1)=0\le d_j(x_j)$, we have $\max_{i\in N} d_i(x_i-1)\le \min_{j\in N} d_j(x_j)$.
It remains to show $\min_{j\in N} d_j(x_j)>0$, because then $\lambda:=\min_{j\in N} d_j(x_j)$ satisfies \eqref{eq:opt-P}.
Assume to the contrary that $d_j(x_j)=0$ for some $j$.
By \eqref{eq:divisor-criterion}, this implies $x_j=0$ and $d(0)=0$.
Then $d_i(x_i-1)\le 0$ for all $i\in N$, which implies $x_i\le 1$ for all $i$, again by \eqref{eq:divisor-criterion}.
Therefore, $\sum_{i\in N}x_i\le n-1<K$, which is a contradiction.
\end{proof}

The fact that every apportionment obtained by the divisor method is optimal for (P) can also be derived from~\cite[Theorem~4.1.1]{IbarakiKatoh1988}, which shows that an optimal solution of a separable convex resource allocation problem is obtained by choosing the $K$ smallest marginal costs.
For (P), choosing the $K$ smallest values among $\{d(k)/p_i \mid i\in N,\ k\in\Z_{\ge 0}\}$ corresponds to assigning seats one by one to a region with the smallest $d(x_i)/p_i$, which is a well-known way of computing divisor methods.

\section{Inverse optimization approach}
\label{sec:inverse}

\subsection{Forward problem: weighted divisor method}

If an observed apportionment deviates from the one obtained by a divisor method, we may interpret that some priority is given to each region.
We represent the relative priority of region $i$ by a nonnegative weight $w_i$, and consider the following problem as the forward problem:
\begin{equation*}
  \text{(FP)} \qquad \mini \ \sum_{i\in N} w_i f_i(x_i) \quad \st \,\, \vt{x}\in\X.
\end{equation*}
We call this method the \emph{weighted divisor method} based on $d$.
The weight $w_i=1$ is the baseline.
If $0\le w_i<1$, the marginal cost of region $i$ decreases, and region $i$ obtains seats more easily.
If $w_i>1$, region $i$ obtains seats less easily.
In fact, for $w_i>0$, the weighted marginal cost of region $i$ is
\begin{equation*}
  w_i\bigl(f_i(x_i+1)-f_i(x_i)\bigr) = w_i d_i(x_i) = \frac{d(x_i)}{p_i/w_i}.
\end{equation*}
Therefore, the weighted divisor method can be interpreted as the divisor method in which the population $p_i$ is replaced by the corrected population $p_i/w_i$.

For $\vt{w}\in\R_{\ge 0}^n$, each $w_i f_i$ is also discrete convex.
Hence, by Lemma~\ref{lem:optimality} and $d_i(-1)=0$, an apportionment $\vt{x}\in\X$ is an optimal solution of (FP) if and only if there exists $\lambda\in\R$ such that
\begin{align}
  w_i d_i(x_i-1) \le \lambda \le w_i d_i(x_i) \qquad (i\in N).
  \label{eq:opt-FP}
\end{align}

\subsection{Formulation of the inverse optimization problem}

Let $\bar{\vt{x}}\in\X$ be an observed apportionment.
Throughout this paper, we make the following assumption.

\begin{assumption}\label{asm:positive}
The observed apportionment satisfies $\bar{x}_i\ge 1$ for all $i\in N$.
\end{assumption}

Assumption~\ref{asm:positive} is natural for many real-world apportionments; for example, every prefecture in Japan has at least one seat.
We assume that $\bar{\vt{x}}$ was obtained by the weighted divisor method based on a given divisor criterion $d$, and we estimate the weight vector $\vt{w}$ from $\bar{\vt{x}}$.
If the resulting optimal weights are far from $\vt{1}$, then $\bar{\vt{x}}$ does not treat the regions equally with respect to the divisor method based on $d$.

Note that the optimality condition \eqref{eq:opt-FP} is invariant under the scaling $(\vt{w},\lambda)\mapsto(c\vt{w},c\lambda)$ for $c>0$.
To remove this degree of freedom, we fix the average of the weights to one, that is, $\sum_{i\in N}w_i=n$.
We then formulate the inverse optimization problem as follows:
\begin{alignat*}{3}
  \text{(IP)} \qquad & \mini \quad && \sum_{i\in N} |w_i-1| \\
  & \st && \sum_{i\in N} w_i = n, \\
  & && d_i(\bar{x}_i-1)\, w_i \le \lambda \le d_i(\bar{x}_i)\, w_i && \quad (i\in N), \\
  & && w_i\ge 0 \quad (i\in N), \qquad \lambda\in\R.
\end{alignat*}
The objective function is the $\ell_1$-distance between $\vt{w}$ and $\vt{1}$, which represents the total amount of correction needed to explain $\bar{\vt{x}}$ by the weighted divisor method.
The second constraint requires that $\bar{\vt{x}}$ is an optimal solution of (FP) under $\vt{w}$ (see \eqref{eq:opt-FP}).

We use the $\ell_1$-norm in the objective function for the following reason.
Our purpose is not only to measure the total deviation from $\vt{1}$, but also to classify each region into two groups: regions that are treated fairly ($w_i=1$) and regions that are not ($w_i\neq 1$).
It is well known that $\ell_1$-norm minimization tends to produce sparse solutions.
In (IP), this sparsity appears in the deviation vector $\vt{w}-\vt{1}$, that is, many components of an optimal $\vt{w}$ are expected to be exactly equal to one.
Thus, for each region $i$, we can judge whether $\bar{\vt{x}}$ treats region $i$ fairly by checking whether $w_i=1$ or not, and the sign of $w_i-1$ shows whether region $i$ is favored ($w_i<1$) or disfavored ($w_i>1$).
In contrast, if we use the squared $\ell_2$-norm $\sum_{i\in N}(w_i-1)^2$, the correction is spread over many regions.
As we will show in Section~\ref{sec:theory}, (IP) indeed has the sparsity property: $w^*_i=1$ holds if and only if $\lambda^*\in[d_i(\bar{x}_i-1),d_i(\bar{x}_i)]$, that is, region $i$ needs no correction at the common parameter $\lambda^*$.
By introducing auxiliary variables $z_i$ $(i\in N)$, replacing the objective function with $\sum_{i\in N}z_i$, and adding the linear constraints $z_i\ge w_i-1$ and $z_i\ge 1-w_i$ $(i\in N)$, (IP) can be written as a linear programming problem.
However, as shown in the next section, (IP) has a special structure that allows a simple and fast algorithm.
We denote the optimal value of (IP) by $\OPT(\bar{\vt{x}})$.

\section{Theoretical properties and algorithm}
\label{sec:theory}

\subsection{Median representation and uniqueness}

For each $i\in N$, let $a_i:=d_i(\bar{x}_i-1),\, b_i:=d_i(\bar{x}_i)$ and $I_i:=[a_i,b_i]$.
By Assumption~\ref{asm:positive} and \eqref{eq:divisor-criterion}, we have $0\le a_i<b_i$.
In particular, $a_i=0$ if and only if $\bar{x}_i=1$ and $d(0)=0$.
Define $\lambda_L:=\max_{i\in N}a_i$ and $\lambda_R:=\min_{i\in N}b_i>0$.
Thus, $\bigcap_{i\in N}I_i$ is nonempty if and only if $\lambda_L\le\lambda_R$.

Every feasible solution of (IP) satisfies $\lambda\ge0$, since $\lambda\ge a_iw_i\ge0$.
We first construct candidate solutions with $\lambda>0$; solutions with $\lambda=0$ remain admissible in (IP).
For a fixed $\lambda>0$, the constraints on $w_i$ are equivalent to
\begin{align*}
  \frac{\lambda}{b_i}\le w_i\le\frac{\lambda}{a_i},
\end{align*}
where $\lambda/a_i:=+\infty$ when $a_i=0$.
The unique minimizer of $|w_i-1|$ over this interval is the projection of $1$ onto the interval, namely
\begin{align}
  w_i(\lambda):=\median\left(
    \frac{\lambda}{b_i},\ 1,\ \frac{\lambda}{a_i}
  \right).
  \label{eq:median}
\end{align}
Equivalently,
\begin{align}
  w_i(\lambda)=
  \begin{cases}
    \lambda/a_i & \text{if }0<\lambda<a_i,\\
    1           & \text{if }a_i\le\lambda\le b_i,\\
    \lambda/b_i & \text{if }\lambda>b_i.
  \end{cases}
  \label{eq:piecewise}
\end{align}
In particular, $w_i(\lambda)<1$ if and only if $\lambda<a_i$, and $w_i(\lambda)>1$ if and only if $\lambda>b_i$.
Let $\vt{w}(\lambda):=(w_1(\lambda),\dots,w_n(\lambda))^\top$ and $S(\lambda):=\sum_{i\in N}w_i(\lambda)$ for any $\lambda>0$.
By construction, $w_i(\lambda)\ge\lambda/b_i>0$ and
$a_iw_i(\lambda)\le\lambda\le b_iw_i(\lambda)$ for every $i\in N$.
Consequently, $(\vt{w}(\lambda),\lambda)$ is feasible for (IP) if and only if $S(\lambda)=n$.
The following lemma describes the structure of this scalar equation; in particular,
it guarantees that $S(\lambda)=n$ always has a positive solution $\lambda^{*}$ and characterizes
when the solution is unique.

\begin{lemma}[Structure of the scalar equation]\label{lem:S}
The following statements hold.
\begin{enumerate}
  \item[(i)] Each $w_i$ is continuous, nondecreasing, and piecewise linear on $(0,+\infty)$, with breakpoints in $\{a_i,b_i\}$.
  Hence, $S$ has the same properties, with breakpoints in $\{a_i,b_i\mid i\in N\}$.
  \item[(ii)] If $\lambda_L\le\lambda_R$, then the set of solutions of $S(\lambda)=n$ is $[\lambda_L,\lambda_R] \cap (0,+\infty)$.
  \item[(iii)] If $\lambda_L>\lambda_R$, then $S$ is strictly increasing on $(0,+\infty)$, and $S(\lambda)=n$ has a unique solution $\lambda^*$ satisfying $\lambda_R<\lambda^*<\lambda_L.$
\end{enumerate}
\end{lemma}

\begin{proof}
(i) This follows directly from \eqref{eq:piecewise}.

(ii) For $\lambda\in[\lambda_L,\lambda_R]\cap(0,+\infty)$, we have $\lambda\in I_i$ for every $i \in N$, so that $w_i(\lambda)=1$ and $S(\lambda)=n$.
If $0<\lambda<\lambda_L$, then $\lambda<b_i$ for every $i \in N$, while $\lambda<a_j$ for at least one $j\in N$.
By \eqref{eq:piecewise}, these imply that $w_i(\lambda)\le 1$ for every $i \in N$ and
$w_j(\lambda)=\lambda/a_j<1$, so that $S(\lambda)<n$.
Similarly, if $\lambda>\lambda_R$, then $\lambda>a_i$ for every $i\in N$ and $\lambda>b_j$ for at least one $j\in N$.
These imply that $w_i(\lambda)\ge 1$ for every $i \in N$ and $w_{j}(\lambda)=\lambda/b_{j}>1$, so that $S(\lambda)>n$.

(iii) Suppose that $S(\lambda_1)=S(\lambda_2)$ for some $0<\lambda_1<\lambda_2$.
Since each $w_i$ is nondecreasing by (i), this implies $w_i(\lambda_1)=w_i(\lambda_2)$ for every $i\in N$.
By \eqref{eq:piecewise}, $w_i$ is strictly increasing on $(0,a_i]$ and on $[b_i,+\infty)$, so such an equality is possible only if $[\lambda_1,\lambda_2]\subseteq I_i$ for every $i\in N$.
This contradicts $\bigcap_{i\in N}I_i=\emptyset$, which is equivalent to $\lambda_L>\lambda_R$.
Hence, $S$ is strictly increasing.

At $\lambda=\lambda_R$, we have $w_i(\lambda_R)\le 1$ for all $i\in N$, since $\lambda_R\le b_i$.
Moreover, for an index $j\in N$ with $a_j=\lambda_L$, we have $0<\lambda_R<\lambda_L=a_j$, and thus $w_j(\lambda_R)=\lambda_R/a_j<1$ by \eqref{eq:piecewise}.
Therefore, $S(\lambda_R)<n$.
Similarly, at $\lambda=\lambda_L$, we have $w_i(\lambda_L)\ge 1$ for all $i\in N$, since $\lambda_L\ge a_i$, and $w_k(\lambda_L)=\lambda_L/b_k>1$ for an index $k\in N$ with $b_k=\lambda_R<\lambda_L$.
Therefore, $S(\lambda_L)>n$.
Since $S$ is continuous by (i) and strictly increasing, the intermediate value theorem yields a unique $\lambda^*\in(\lambda_R,\lambda_L)$ with $S(\lambda^*)=n$.
\end{proof}

For $\vt{w}\in\R_{\ge0}^n$, define
$L(\vt{w}):=\{i\in N\mid w_i<1\}$ and $U(\vt{w}):=\{i\in N\mid w_i>1\}$.
If $\sum_{i\in N}(w_i-1)=0$, then
\begin{align}
  \delta(\vt{w}):=\sum_{i\in U(\vt{w})}(w_i-1)=\sum_{i\in L(\vt{w})}(1-w_i)
  \quad\text{and}\quad
  \sum_{i\in N}|w_i-1|=2\delta(\vt{w}).
  \label{eq:split}
\end{align}

The next theorem shows that any positive solution $\lambda^*$ of $S(\lambda)=n$ yields an optimal solution of (IP); that is, enforcing the normalization $S(\lambda)=n$ on the componentwise
projections solves (IP).

\begin{theorem}[Median representation and uniqueness]\label{thm:median}
Let $\lambda^*>0$ satisfy $S(\lambda^*)=n$, and let $\vt{w}^*:=\vt{w}(\lambda^*)$.
Then $(\vt{w}^*,\lambda^*)$ is an optimal solution of (IP), and $\vt{w}^*$ is the unique optimal weight vector.
Moreover, $w_i^*>0$ for every $i\in N$.
\end{theorem}

\begin{proof}
Since $(\vt{w}^*,\lambda^*)$ is feasible by construction, \eqref{eq:split} holds for $\vt{w}^*$.
If $\delta(\vt{w}^*)=0$, then $\vt{w}^*=\vt{1}$, which is clearly the unique weight vector
attaining the minimum objective value of zero.
We therefore assume $\delta(\vt{w}^*)>0$, so that both $L(\vt{w}^*)$ and $U(\vt{w}^*)$ are nonempty.

Let $(\vt{w},\lambda)$ be any feasible solution of (IP), allowing $\lambda=0$.
If $\lambda\ge\lambda^*$, then for each $i\in U(\vt{w}^*)$,
\begin{align}
  w_i\ge\frac{\lambda}{b_i}\ge\frac{\lambda^*}{b_i}=w_i^*\ (>1)
  \label{eq:proof-U}
\end{align}
by the constraint $\lambda\le b_iw_i$ and \eqref{eq:piecewise}.
In particular, $U(\vt{w}^*)\subseteq U(\vt{w})$, and hence, by \eqref{eq:split},
\begin{align*}
  \sum_{i\in N}|w_i-1|=2\sum_{i\in U(\vt{w})}(w_i-1)\ge2\sum_{i\in U(\vt{w}^*)}(w_i-1)
  \ge2\sum_{i\in U(\vt{w}^*)}(w_i^*-1)=\sum_{i\in N}|w_i^*-1|.
\end{align*}
If $\lambda\le\lambda^*$, then for each $i\in L(\vt{w}^*)$, we have $a_i>\lambda^*>0$
by \eqref{eq:piecewise}, and the constraint $a_iw_i\le\lambda$ gives
\begin{align}
  w_i\le\frac{\lambda}{a_i}\le\frac{\lambda^*}{a_i}=w_i^*\ (<1),
  \label{eq:proof-L}
\end{align}
which is valid also for $\lambda=0$.
Since $L(\vt{w}^*)\subseteq L(\vt{w})$, we have, by \eqref{eq:split},
\begin{equation*}
  \sum_{i\in N}|w_i-1|=2\sum_{i\in L(\vt{w})}(1-w_i)\ge2\sum_{i\in L(\vt{w}^*)}(1-w_i)
  \ge2\sum_{i\in L(\vt{w}^*)}(1-w_i^*)=\sum_{i\in N}|w_i^*-1|.
\end{equation*}
These inequalities prove the optimality of $(\vt{w}^*,\lambda^*)$.
Moreover, if $\lambda>\lambda^*$ (resp.\ $\lambda<\lambda^*$), then the second inequality
in \eqref{eq:proof-U} (resp.\ \eqref{eq:proof-L}) is strict for every $i \in U(\vt{w}^*)$ (resp.\ $L(\vt{w}^*)$), and thus $\sum_{i\in N}|w_i-1|>\sum_{i\in N}|w_i^*-1|$.
Therefore, every optimal solution satisfies $\lambda=\lambda^*$.
For this $\lambda^*$, each $w_i^*$ is the unique minimizer of $|w_i-1|$ over the feasible interval $[\lambda^*/b_i,\lambda^*/a_i]$, which proves the uniqueness of $\vt{w}^*$.
\end{proof}

Theorem~\ref{thm:median} reduces (IP) to the scalar equation $S(\lambda)=n$.
This is analogous to the forward divisor method, where a divisor is chosen to make the total number of allocated seats equal to $K$.
The proof also shows that $\OPT(\bar{\vt{x}})=2\delta(\vt{w}^*)$.
The following corollary characterizes when no correction is needed.

\begin{corollary}\label{cor:lambda}
The following conditions are equivalent:
(a) $\OPT(\bar{\vt{x}})=0$, that is, the optimal weight vector is $\vt{w}^*=\vt{1}$;
(b) $\lambda_L\le\lambda_R$, or equivalently $\bigcap_{i\in N}I_i\ne\emptyset$;
(c) $\bar{\vt{x}}$ is obtained by the unweighted divisor method based on $d$.
\end{corollary}

\begin{proof}
The objective value is zero exactly when $\vt{w}=\vt{1}$, and $(\vt{1},\lambda)$ is feasible exactly when $\lambda\in\bigcap_{i\in N}I_i$.
This proves the equivalence of (a) and (b).
Since $\lambda_R>0$, a nonempty intersection contains a positive number.
The equivalence of (b) and (c) therefore follows from Theorem~\ref{thm:divisor-P}.
\end{proof}

\begin{example}\label{ex:small}
We illustrate Theorem~\ref{thm:median} by a small example.
Consider three regions A, B, and C with populations $p=(610,350,40)$ and $K=10$, and let $d$ be the divisor criterion of the Adams method, $d(x)=x$.
The Adams method gives the apportionment $(6,3,1)$; for example, the divisor $r=120$ yields $\lceil 610/120\rceil+\lceil 350/120\rceil+\lceil 40/120\rceil=6+3+1=10$.
Suppose that the apportionment $\bar{\vt{x}}=(5,4,1)$ is observed, that is, one seat of A is moved to B.
Since $d(x)=x$, we have $a_i=(\bar{x}_i-1)/p_i$ and $b_i=\bar{x}_i/p_i$, and thus
\begin{equation*}
  I_A=\Bigl[\frac{4}{610},\,\frac{5}{610}\Bigr],\qquad
  I_B=\Bigl[\frac{3}{350},\,\frac{4}{350}\Bigr],\qquad
  I_C=\Bigl[0,\,\frac{1}{40}\Bigr].
\end{equation*}
The interval $I_A$ lies entirely to the left of $I_B$, that is, $\lambda_L=3/350>1/122=\lambda_R$.
Hence, $\bigcap_{i\in N}I_i=\emptyset$, and $\bar{\vt{x}}$ is not obtained by the Adams method (Corollary~\ref{cor:lambda}).
For $\lambda\in[1/122,\,3/350]$, we have $w_A(\lambda)=122\lambda$, $w_B(\lambda)=(350/3)\lambda$, and $w_C(\lambda)=1$ by \eqref{eq:piecewise}.
Thus, the equation $S(\lambda)=(716/3)\lambda+1=3$ gives $\lambda^*=3/358\approx 0.00838$, which is in this interval.
The unique optimal weight vector is
\begin{equation*}
  \vt{w}^*=\left(\frac{183}{179},\ \frac{175}{179},\ 1\right)\approx(1.022,\ 0.978,\ 1),
\end{equation*}
and $\OPT(\bar{\vt{x}})=8/179\approx 0.045$.
The weight of B, which receives one more seat than in the Adams apportionment, is less than one, and the weight of A, which receives one less seat, is greater than one.
On the other hand, the weight of C is exactly one.
This shows that, by the sparsity of the $\ell_1$-norm, the region C is judged to be treated fairly, and the correction is concentrated on the regions A and B.
We also note that, under $\vt{w}^*$, the Adams apportionment $(6,3,1)$ is also optimal for (FP), since $w^*_B a_B=\lambda^*=w^*_A b_A$ and the exchange of one seat from B to A does not change the objective value.
In general, if $\OPT(\bar{\vt{x}})>0$, then for any $i,j\in N$ with $w^*_i<1<w^*_j$, we have $w^*_i a_i=\lambda^*=w^*_j b_j$ by \eqref{eq:piecewise}, and thus $\bar{\vt{x}}-\vt{e}_i+\vt{e}_j$ is also optimal for (FP) under $\vt{w}^*$.
That is, $\bar{\vt{x}}$ is not the unique optimal solution of (FP) under $\vt{w}^*$.
\end{example}

\subsection{Exchange property of the optimal weights}

In this subsection, we show that the optimal weights of (IP) indicate a direction in which an apportionment should be corrected.
Let $\X_+ = \{\vt{x}\in\X \mid x_i\ge 1\ (i\in N)\}$,
that is, the set of apportionments satisfying the condition of Assumption~\ref{asm:positive}. We consider (IP) for each $\vt{x}\in\X_+$ in place of $\bar{\vt{x}}$.
We denote by $\OPT(\vt{x})$ the optimal value and by $\vt{w}^*(\vt{x})$ the unique optimal weight vector (Theorem~\ref{thm:median}).
We define
\begin{equation*}
  N^-(\vt{x}):=\{i\in N\mid w^*_i(\vt{x})<1\},\qquad
  N^+(\vt{x}):=\{i\in N\mid w^*_i(\vt{x})>1\},
\end{equation*}
that is, $N^-(\vt{x})=L(\vt{w}^*(\vt{x}))$ and $N^+(\vt{x})=U(\vt{w}^*(\vt{x}))$.
Since $\sum_{i\in N}w^*_i(\vt{x})=n$, both $N^-(\vt{x})$ and $N^+(\vt{x})$ are nonempty if $\OPT(\vt{x})>0$.
We also denote the objective function of (P), that is, the unweighted divisor method, by $F(\vt{x}):=\sum_{i\in N}f_i(x_i)$.

\begin{theorem}\label{thm:exchange}
Let $\vt{x}\in\X_+$ with $\OPT(\vt{x})>0$, and let $i\in N^-(\vt{x})$ and $j\in N^+(\vt{x})$.
Assume that $x_i\ge 2$, and let $\vt{x}':=\vt{x}-\vt{e}_i+\vt{e}_j$.
Then $\vt{x}'\in\X_+$, and
\begin{equation*}
  \OPT(\vt{x}')<\OPT(\vt{x}) \qquad\text{and}\qquad F(\vt{x}')<F(\vt{x}).
\end{equation*}
\end{theorem}

\begin{proof}
Let $\vt{w}^*:=\vt{w}^*(\vt{x})$, and let $\lambda^*>0$ satisfy $S(\lambda^*)=n$ for $\vt{x}$.
Since $x_i\ge 2$, we have $\vt{x}'\in\X_+$.
By \eqref{eq:piecewise}, $i\in N^-(\vt{x})$ and $j\in N^+(\vt{x})$ imply
\begin{align}
  w^*_i=\frac{\lambda^*}{d_i(x_i-1)}<1, \qquad w^*_j=\frac{\lambda^*}{d_j(x_j)}>1,
  \label{eq:exchange-LU}
\end{align}
and hence $d_i(x_i-1)>\lambda^*>d_j(x_j)$.

\emph{Decrease of $F$.}
Since $f_k(x_k+1)-f_k(x_k)=d_k(x_k)$, we have $F(\vt{x}')-F(\vt{x}) = d_j(x_j)-d_i(x_i-1) < 0$.

\emph{Decrease of $\OPT$.}
Let $0<\varepsilon<\min\{1-w^*_i,\ w^*_j-1\}$, and define $\tilde{\vt{w}}:=\vt{w}^*+\varepsilon\vt{e}_i-\varepsilon\vt{e}_j$.
We show that $(\tilde{\vt{w}},\lambda^*)$ is feasible for (IP) with respect to $\vt{x}'$ if $\varepsilon$ is sufficiently small.
Clearly $\sum_{k\in N}\tilde{w}_k=n$ and $\tilde{\vt{w}}\ge\vt{0}$.
For $k\neq i,j$, the constraints do not change.
For region $i$, the constraint is $d_i(x_i-2)\,\tilde{w}_i\le\lambda^*\le d_i(x_i-1)\,\tilde{w}_i$.
By \eqref{eq:exchange-LU}, $d_i(x_i-1)\,\tilde{w}_i=\lambda^*+\varepsilon\, d_i(x_i-1)\ge\lambda^*$.
Since $d$ is strictly increasing, $d_i(x_i-2)\,w^*_i<d_i(x_i-1)\,w^*_i=\lambda^*$, and thus $d_i(x_i-2)\,\tilde{w}_i\le\lambda^*$ for sufficiently small $\varepsilon$.
Similarly, for region $j$, the constraint is $d_j(x_j)\,\tilde{w}_j\le\lambda^*\le d_j(x_j+1)\,\tilde{w}_j$.
We have $d_j(x_j)\,\tilde{w}_j=\lambda^*-\varepsilon\, d_j(x_j)\le\lambda^*$, and $d_j(x_j+1)\,\tilde{w}_j\ge\lambda^*$ for sufficiently small $\varepsilon$, since $d_j(x_j+1)\,w^*_j>d_j(x_j)\,w^*_j=\lambda^*$.
Therefore, $(\tilde{\vt{w}},\lambda^*)$ is feasible.
Since $w^*_i<\tilde{w}_i<1<\tilde{w}_j<w^*_j$, we obtain $\OPT(\vt{x}') \le \sum_{k\in N}|\tilde{w}_k-1| = \OPT(\vt{x})-2\varepsilon < \OPT(\vt{x})$.
\end{proof}

Theorem~\ref{thm:exchange} shows that moving one seat from a favored region ($w^*_i<1$) to a disfavored region ($w^*_j>1$) improves the fit to the chosen divisor method, in terms of both the inverse optimization problem and the original problem (P).
In Example~\ref{ex:small}, we have $N^-(\bar{\vt{x}})=\{\mathrm{B}\}$ and $N^+(\bar{\vt{x}})=\{\mathrm{A}\}$, and the exchange gives the Adams apportionment $(6,3,1)$ with $\OPT=0$.

The assumption $x_i\ge 2$ is needed to keep $\vt{x}'\in\X_+$.
If $d(0)=0$, as in the Adams, Dean, and Hill methods, this assumption always holds.
Indeed, $i\in N^-(\vt{x})$ implies $d_i(x_i-1)>\lambda^*>0$ by \eqref{eq:exchange-LU}, and hence $x_i\ge 2$.
As a natural consequence of Theorem~\ref{thm:exchange}, if $d(0)=0$, we can repeat the exchange: starting from any $\vt{x}\in\X_+$, we move one seat from a region in $N^-(\vt{x})$ to a region in $N^+(\vt{x})$ and recompute the optimal weights, as long as $\OPT(\vt{x})>0$.
Since $F$ strictly decreases at each step, no apportionment appears twice, and since $\X_+$ is finite, this procedure terminates after a finite number of exchanges.
It stops only when $\OPT(\vt{x})=0$, because $N^-(\vt{x})$ and $N^+(\vt{x})$ are nonempty if $\OPT(\vt{x})>0$; that is, the final apportionment is obtained by the divisor method based on $d$ (Corollary~\ref{cor:lambda}).

\begin{remark}\label{rem:Mconvex}
As noted in Section~\ref{sec:prelim}, the function $F$ restricted to $\X$ is an M-convex function. Hence, the above procedure can be regarded as an exchange-based descent method for minimizing $F$.
\end{remark}

\subsection{Algorithm}

By Lemma~\ref{lem:S}, $S(\lambda)$ is a nondecreasing piecewise linear function whose breakpoints are in $\{a_i,b_i\mid i\in N\}$.
Therefore, we can find $\lambda^*$ by sorting the breakpoints, finding the interval containing $\lambda^*$ by binary search, and solving a linear equation in that interval.
The algorithm is described in Algorithm~\ref{alg:IP}.

\begin{algorithm}[hbt]
\caption{Algorithm for (IP)}
\label{alg:IP}
\begin{algorithmic}[1]
\Require Populations $p_i$, observed apportionment $\bar{\vt{x}}$ satisfying Assumption~\ref{asm:positive}, divisor criterion $d$
\Ensure Optimal weight vector $\vt{w}^*$ of (IP)
\State Compute $a_i\gets d(\bar{x}_i-1)/p_i$ and $b_i\gets d(\bar{x}_i)/p_i$ for each $i\in N$
\State $\lambda_L\gets\max_{i\in N}a_i$, \ $\lambda_R\gets\min_{i\in N}b_i$
\If{$\lambda_L\le\lambda_R$}
  \State \Return $\vt{w}^*=\vt{1}$ \Comment{$\OPT(\bar{\vt{x}})=0$ by Corollary~\ref{cor:lambda}}
\EndIf
\State Sort the distinct positive values in $\{a_i,b_i\mid i\in N\}$ as $s_1<s_2<\dots<s_m$, and set $s_0\gets 0$
\State Find an index $\ell\in\{0,1,\dots,m-1\}$ with $S(s_\ell)\le n\le S(s_{\ell+1})$ by binary search, where $S(s_0):=|\{i\in N\mid a_i=0\}|$
\State On $[s_\ell,s_{\ell+1}]$, write $S(\lambda)=\alpha\lambda+\beta$ with
\Statex \qquad $\displaystyle \alpha\gets\sum_{i:\,a_i\ge s_{\ell+1}}\frac{1}{a_i}+\sum_{i:\,b_i\le s_\ell}\frac{1}{b_i}, \qquad \beta\gets|\{i\in N\mid a_i\le s_\ell,\ b_i\ge s_{\ell+1}\}|$
\State $\lambda^*\gets(n-\beta)/\alpha$
\State \Return $\vt{w}^*=(w_1(\lambda^*),\dots,w_n(\lambda^*))^\top$ computed by \eqref{eq:median}
\end{algorithmic}
\end{algorithm}

\begin{theorem}\label{thm:complexity}
Algorithm~\ref{alg:IP} computes the unique optimal weight vector $\vt{w}^*$ of (IP) in $O(n\log n)$ time.
\end{theorem}

\begin{proof}
If $\lambda_L\le\lambda_R$, the output is correct by Corollary~\ref{cor:lambda}.
Otherwise, $S$ is strictly increasing by Lemma~\ref{lem:S}\,(iii).
By \eqref{eq:piecewise} and $\lambda_L>\lambda_R>0$, we have $\lim_{\lambda\to+0}S(\lambda)=S(s_0)<n$ and $S(s_m)\ge n$, where $s_m=\max_{i\in N}b_i$.
Hence, an index $\ell$ in Line~7 exists, and $\lambda^*\in[s_\ell,s_{\ell+1}]$.
Since there is no breakpoint in the open interval $(s_\ell,s_{\ell+1})$, \eqref{eq:piecewise} shows that $S(\lambda)=\alpha\lambda+\beta$ holds on $[s_\ell,s_{\ell+1}]$ with $\alpha$ and $\beta$ in Line~8.
Since $S$ is strictly increasing, $\alpha>0$, and Line~9 gives the solution $\lambda^*$ of $S(\lambda)=n$.
By Theorem~\ref{thm:median}, $\vt{w}(\lambda^*)$ is the unique optimal weight vector.

Lines~1--5, 8, 9, and 10 take $O(n)$ time.
Sorting in Line~6 takes $O(n\log n)$ time.
In Line~7, each evaluation of $S$ takes $O(n)$ time, and the binary search requires $O(\log n)$ evaluations, so Line~7 takes $O(n\log n)$ time.
Therefore, the total running time is $O(n\log n)$.
\end{proof}

\begin{remark}
Since $S(\lambda)$ is the sum of $n$ nondecreasing piecewise linear functions with at most two breakpoints each, the running time can be reduced to $O(n)$ by the breakpoint searching technique of Brucker~\cite{Brucker1984} for the continuous quadratic knapsack problem, which uses linear-time selection~\cite{BlumEtAl1973}.
We adopt the simpler algorithm based on sorting, since $n$ is small in typical applications (for example, $n=47$ prefectures for the House of Representatives of Japan and $n=50$ states for the House of Representatives of the United States).
\end{remark}

\section{Application to the House of Representatives of Japan}
\label{sec:experiment}

We apply the proposed approach to the apportionment of the seats of the single-member districts of the House of Representatives of Japan among the $n=47$ prefectures.

\paragraph{Background and data}
The current electoral system was introduced in 1994, and the seats of the single-member districts were apportioned by the one-seat-per-prefecture rule with $K=300$.
Since then, the apportionment has been revised several times: the ``5 increase 5 decrease'' revision in 2002, the ``0 increase 5 decrease'' revision in 2013 ($K=295$), the ``0 increase 6 decrease'' revision in 2017 ($K=289$), and the ``10 increase 10 decrease'' revision in 2022 ($K=289$).
The last one is the first apportionment based on the Adams method.
We consider the seven general elections held in 2005, 2009, 2012, 2014, 2017, 2021, and 2024, and we use the number of seats of each prefecture in election year $t$ as the observed apportionment $\bar{\vt{x}}^t$.
The total number of seats is $K=300$ for 2005--2012, $K=295$ for 2014, and $K=289$ for 2017--2024.
As the population of each prefecture, we use the most recent census before each election: the total population for the elections up to 2014, and the Japanese population for the elections from 2017 onward.
Hence, the same population data are used for 2012 and 2014, and for 2021 and 2024.
For each of the five divisor methods $m$ in Section~\ref{sec:prelim} and each election year $t$, we solved (IP) by Algorithm~\ref{alg:IP} and computed the optimal value $\OPT_{m,t}:=\sum_{i\in N}|w^*_{i,m,t}-1|$.
The computation time was negligible in all cases.
The optimal weights of all prefectures for every election and method are listed
in~\ref{app:weights}.

\begin{table}[t]
  \centering
  \caption{Optimal values $\OPT_{m,t}$ of (IP) for each divisor method $m$ and election year $t$.}
  \label{tab:opt}
  \begin{tabular}{lrrrrrrr}
    \toprule
    Method & 2005 & 2009 & 2012 & 2014 & 2017 & 2021 & 2024 \\
    \midrule
    Adams     & 0.703 & 0.877 & 1.189 & 0.529 & 0.355 & 0.614 & 0.000 \\
    Dean      & 1.790 & 1.991 & 2.332 & 1.447 & 1.075 & 1.427 & 0.358 \\
    Hill      & 1.923 & 2.129 & 2.493 & 1.591 & 1.228 & 1.578 & 0.504 \\
    Webster   & 2.057 & 2.293 & 2.695 & 1.752 & 1.391 & 1.742 & 0.665 \\
    Jefferson & 3.564 & 3.890 & 4.415 & 3.261 & 2.940 & 3.422 & 2.244 \\
    \bottomrule
  \end{tabular}
\end{table}

\paragraph{Comparison of divisor methods}
Table~\ref{tab:opt} shows the optimal values.
In all the seven elections, the optimal values satisfy
\begin{equation*}
  \OPT_{\mathrm{Adams},t}<\OPT_{\mathrm{Dean},t}<\OPT_{\mathrm{Hill},t}<\OPT_{\mathrm{Webster},t}<\OPT_{\mathrm{Jefferson},t}.
\end{equation*}
That is, a divisor method that is more favorable to less populated prefectures explains the observed apportionments with a smaller amount of correction.
This is consistent with the fact that the one-seat-per-prefecture rule favors less populated prefectures.
In 2024, $\OPT_{\mathrm{Adams},t}=0$ holds, since the actual apportionment is based on the Adams method (Corollary~\ref{cor:lambda}).
We also observe the following tendency for the Adams method.
When the population data are updated and the apportionment is not changed, the optimal value increases (from 0.703 to 1.189 in 2005--2012, and from 0.355 to 0.614 in 2017--2021).
When the apportionment is revised and the population data are not changed, the optimal value decreases (from 1.189 to 0.529 in 2012--2014, and from 0.614 to 0 in 2021--2024).

\paragraph{Optimal weights}
The optimal weights are sparse, as intended by the $\ell_1$ objective
(Section~\ref{sec:inverse}): for the Adams method, 29 to 34 of the 47 prefectures
have $w_i^*=1$ in each election from 2005 to 2021, and the correction concentrates
on the remaining prefectures (\ref{app:weights}).
For the Adams method, metropolitan prefectures have weights greater than one, and less populated prefectures have weights less than one.
For example, in 2012, the five largest weights are those of Tokyo (1.161), Kanagawa (1.109), Aichi (1.090), Saitama (1.058), and Chiba (1.055), and the five smallest weights are those of Kochi (0.843), Tokushima (0.866), Fukui (0.889), Saga (0.937), and Kagoshima (0.941).
A similar tendency is observed for the other four methods (see~\ref{app:weights}).

\paragraph{Consistency with the revisions}
We compare the optimal weights with the revisions of the apportionment made in 2013, 2017, and 2022.
In the ``0 increase 5 decrease'' revision, the number of seats was reduced in Yamanashi, Fukui, Tokushima, Kochi, and Saga.
In the ``0 increase 6 decrease'' revision, it was reduced in Aomori, Iwate, Mie, Nara, Kumamoto, and Kagoshima.
For all the five methods, all these 11 prefectures have $w^*_i<1$ in the apportionment just before the revision (2012 and 2014, respectively).
Although these revisions change $K$ and are not covered by Theorem~\ref{thm:exchange}, the seats were reduced only in prefectures that are favored according to the optimal weights.

The ``10 increase 10 decrease'' revision keeps $K=289$, and it can be regarded as a set of one-seat exchanges.
The number of seats was increased in Tokyo ($+5$), Kanagawa ($+2$), Saitama, Chiba, and Aichi ($+1$ each), and decreased by one in Miyagi, Fukushima, Niigata, Shiga, Wakayama, Okayama, Hiroshima, Yamaguchi, Ehime, and Nagasaki.
In the 2021 apportionment, for all the five methods, the five prefectures with increased seats have $w^*_i>1$, and the ten prefectures with decreased seats have $w^*_i\le 1$.
For the Adams method, all the ten prefectures have $w^*_i<1$; that is, all the exchanges are from $N^-(\bar{\vt{x}}^{2021})$ to $N^+(\bar{\vt{x}}^{2021})$, which is the type of exchange considered in Theorem~\ref{thm:exchange}.
For the other methods, a few prefectures have $w^*_i=1$ (Hiroshima for the Dean, Hill, and Webster methods, and Miyagi and Hiroshima for the Jefferson method; see Table~\ref{tab:w2021}).
Therefore, no seat was moved in the direction opposite to that indicated by the optimal weights, and the optimal value decreased for all the five methods (Table~\ref{tab:opt}).
These results show that the optimal weights can explain the direction of the three revisions.

\section{Conclusion}
\label{sec:conclusion}

In this paper, we proposed an inverse optimization approach to evaluate an observed apportionment with respect to a divisor method.
We introduced the weighted divisor method as the forward problem, and formulated the inverse optimization problem (IP) that finds the weight vector closest to $\vt{1}$ in the $\ell_1$-norm under which the observed apportionment is optimal.
For a general divisor method, we showed that the optimal weight vector is unique and is expressed as the median of three values that depend on a single parameter $\lambda$.
By this result, (IP) is reduced to the one-dimensional monotone equation $S(\lambda)=n$, and we gave an $O(n\log n)$ time algorithm.
We also showed that moving one seat from a region with $w^*_i<1$ to a region with $w^*_j>1$ strictly decreases both the optimal value of (IP) and the objective value of the unweighted divisor method.
Furthermore, we applied the proposed approach to the apportionment of the House of Representatives of Japan.
The results show that divisor methods more favorable to less populated prefectures explain the past apportionments with smaller corrections, and that the revisions actually made are consistent with the directions indicated by the optimal weights.

A direction for future research is to extend the model to apportionment problems with additional constraints, such as lower and upper bounds on the number of seats of each region (e.g., in the European Parliament) or on the total number of seats of groups of regions.
If these constraints are defined on a laminar family, the forward problem is still a separable discrete convex resource allocation problem~\cite{Hochbaum1994,KatohShiouraIbaraki2013}, but its optimality condition cannot be described by a single parameter $\lambda$ in general.
It is an important question whether Theorem~\ref{thm:median} can be extended to such problems.

\section*{Acknowledgements}
This work was supported by JSPS KAKENHI Grant Numbers JP22K17856 and JP26K14719.

\section*{Declaration of competing interest}
The authors declare that they have no known competing financial interests or personal relationships that could have appeared to influence the work reported in this paper.

\section*{Declaration of generative AI and AI-assisted technologies in the writing process}
During the preparation of this work, the authors used Claude Opus~5.5 (Anthropic) in order to improve the language and readability of the manuscript.
After using this tool, the authors reviewed and edited the content as needed and take full responsibility for the content of the publication.

\section*{Data availability}
The census data and the numbers of seats used in this study are publicly available from official sources.
The optimal weights computed in this study are provided in the supplementary material.

\bibliographystyle{abbrvnat}
\bibliography{refs}

\begin{thebibliography}{12}
\providecommand{\natexlab}[1]{#1}
\providecommand{\url}[1]{\texttt{#1}}
\expandafter\ifx\csname urlstyle\endcsname\relax
  \providecommand{\doi}[1]{doi: #1}\else
  \providecommand{\doi}{doi: \begingroup \urlstyle{rm}\Url}\fi

\bibitem[Ahuja and Orlin(2001)]{AhujaOrlin2001}
R.~K. Ahuja and J.~B. Orlin.
\newblock Inverse optimization.
\newblock \emph{Operations Research}, 49:\penalty0 771--783, 2001.

\bibitem[Balinski and Young(2010)]{BalinskiYoung2001}
M.~L. Balinski and H.~P. Young.
\newblock \emph{Fair Representation: Meeting the Ideal of {One Man, One Vote}}.
\newblock Bloomsbury Publishing, second edition, 2010.

\bibitem[Blum et~al.(1973)Blum, Floyd, Pratt, Rivest, and Tarjan]{BlumEtAl1973}
M.~Blum, R.~W. Floyd, V.~Pratt, R.~L. Rivest, and R.~E. Tarjan.
\newblock Time bounds for selection.
\newblock \emph{Journal of Computer and System Sciences}, 7:\penalty0 448--461, 1973.

\bibitem[Brucker(1984)]{Brucker1984}
P.~Brucker.
\newblock An {$O(n)$} algorithm for quadratic knapsack problems.
\newblock \emph{Operations Research Letters}, 3:\penalty0 163--166, 1984.

\bibitem[Chan et~al.(2025)Chan, Mahmood, and Zhu]{ChanMahmoodZhu2025}
T.~C.~Y. Chan, R.~Mahmood, and I.~Y. Zhu.
\newblock Inverse optimization: Theory and applications.
\newblock \emph{Operations Research}, 73:\penalty0 1046--1074, 2025.

\bibitem[Heuberger(2004)]{Heuberger2004}
C.~Heuberger.
\newblock Inverse combinatorial optimization: A survey on problems, methods, and results.
\newblock \emph{Journal of Combinatorial Optimization}, 8:\penalty0 329--361, 2004.

\bibitem[Hochbaum(1994)]{Hochbaum1994}
D.~S. Hochbaum.
\newblock Lower and upper bounds for the allocation problem and other nonlinear optimization problems.
\newblock \emph{Mathematics of Operations Research}, 19:\penalty0 390--409, 1994.

\bibitem[Ibaraki and Katoh(1988)]{IbarakiKatoh1988}
T.~Ibaraki and N.~Katoh.
\newblock \emph{Resource Allocation Problems: Algorithmic Approaches}.
\newblock MIT Press, Cambridge, MA, 1988.

\bibitem[Katoh et~al.(2013)Katoh, Shioura, and Ibaraki]{KatohShiouraIbaraki2013}
N.~Katoh, A.~Shioura, and T.~Ibaraki.
\newblock Resource allocation problems.
\newblock In P.~M. Pardalos, D.-Z. Du, and R.~L. Graham, editors, \emph{Handbook of Combinatorial Optimization}, pages 2897--2988. Springer, New York, second edition, 2013.

\bibitem[Murota(2003)]{Murota2003}
K.~Murota.
\newblock \emph{Discrete Convex Analysis}.
\newblock SIAM, Philadelphia, 2003.

\bibitem[Murota(2024)]{Murota2024}
K.~Murota.
\newblock \emph{Discrete Convex Analysis: Extended Theory with Applications}.
\newblock Maruzen Publishing, Tokyo, 2024.
\newblock (in Japanese).

\bibitem[Pukelsheim(2017)]{Pukelsheim2017}
F.~Pukelsheim.
\newblock \emph{Proportional Representation: Apportionment Methods and Their Applications}.
\newblock Springer, Cham, second edition, 2017.

\end{thebibliography}

\newpage
%

\appendix
\section{Supplementary Material: Optimal weights of the prefectures}\label{app:weights}

This appendix lists the optimal weights $w_i^*$ of all prefectures for each election
year and each divisor method, computed in Section~5.
Tables~\ref{tab:w2005}--\ref{tab:w2024} share the following conventions.
\begin{itemize}
  \item $p_i$ and $\bar{x}_i$ denote the population and the number of seats of
  prefecture $i$, respectively.
  \item A value of exactly one is written as 1; other values are rounded to three
  decimal places. Thus, 1.000 means that $w_i^*\ne1$.
  \item The superscripts $+$ and $-$ attached to a prefecture indicate an increase
  and a decrease in its number of seats in the next election, respectively.
  \item The population is the total population in the census for the elections up
  to 2014, and the Japanese population in the census for the elections from 2017
  onward (Table~\ref{tab:data}).
\end{itemize}

\begin{table}[h]
  \centering
  \caption{Total number of seats and population data for each election.}
  \label{tab:data}
  \begin{tabular}{lll}
    \toprule
    Election & $K$ & Population \\
    \midrule
    2005       & 300 & 2000 census, total population \\
    2009       & 300 & 2005 census, total population \\
    2012       & 300 & 2010 census, total population \\
    2014       & 295 & 2010 census, total population \\
    2017       & 289 & 2015 census, Japanese population \\
    2021, 2024 & 289 & 2020 census, Japanese population \\
    \bottomrule
  \end{tabular}
\end{table}

\begin{table}[htbp]
  \centering
  \footnotesize
  \caption{Optimal weights $w_i^*$ of the prefectures for the 2005 election
($K=300$).}
  \label{tab:w2005}
  \begin{tabular}{lrrrrrrr}
    \toprule
    \textbf{Prefecture} & \textbf{$p_i$} & \textbf{$\bar{x}_i$} & \textbf{Adams} & \textbf{Dean} & \textbf{Hill} & \textbf{Webster} & \textbf{Jefferson} \\
    \midrule
    Hokkaido & 5,683,062 & 12 & 1.050 & 1.100 & 1.107 & 1.114 & 1.181 \\
    Aomori & 1,475,728 & 4 & 1 & 1 & 1 & 1 & 0.997 \\
    Iwate & 1,416,180 & 4 & 1 & 0.998 & 0.995 & 0.991 & 0.957 \\
    Miyagi & 2,365,320 & 6 & 1 & 1 & 1 & 1 & 1 \\
    Akita & 1,189,279 & 3 & 1 & 1 & 1 & 1 & 1 \\
    Yamagata & 1,244,147 & 3 & 1 & 1 & 1 & 1 & 1 \\
    Fukushima & 2,126,935 & 5 & 1 & 1 & 1 & 1 & 1 \\
    Ibaraki & 2,985,676 & 7 & 1 & 1 & 1 & 1 & 1.009 \\
    Tochigi & 2,004,817 & 5 & 1 & 1 & 1 & 1 & 1 \\
    Gunma & 2,024,852 & 5 & 1 & 1 & 1 & 1 & 1 \\
    Saitama & 6,938,006 & 15 & 1.026 & 1.083 & 1.090 & 1.096 & 1.172 \\
    Chiba & 5,926,285 & 13 & 1.011 & 1.062 & 1.069 & 1.075 & 1.144 \\
    Tokyo & 12,064,101 & 25 & 1.070 & 1.144 & 1.151 & 1.159 & 1.254 \\
    Kanagawa & 8,489,974 & 18 & 1.046 & 1.110 & 1.117 & 1.124 & 1.208 \\
    Niigata & 2,475,733 & 6 & 1 & 1 & 1 & 1 & 1 \\
    Toyama & 1,120,851 & 3 & 1 & 1 & 1 & 1 & 1 \\
    Ishikawa & 1,180,977 & 3 & 1 & 1 & 1 & 1 & 1 \\
    Fukui & 828,944 & 3 & 0.919 & 0.835 & 0.823 & 0.812 & 0.747 \\
    Yamanashi & 888,172 & 3 & 0.985 & 0.894 & 0.882 & 0.870 & 0.800 \\
    Nagano & 2,215,168 & 5 & 1 & 1 & 1 & 1 & 1 \\
    Gifu & 2,107,700 & 5 & 1 & 1 & 1 & 1 & 1 \\
    Shizuoka & 3,767,393 & 8 & 1.044 & 1.075 & 1.080 & 1.086 & 1.131 \\
    Aichi & 7,043,300 & 15 & 1.041 & 1.099 & 1.106 & 1.113 & 1.190 \\
    Mie & 1,857,339 & 5 & 1 & 1 & 1 & 1 & 1 \\
    Shiga & 1,342,832 & 4 & 0.992 & 0.946 & 0.943 & 0.940 & 0.907 \\
    Kyoto & 2,644,391 & 6 & 1 & 1 & 1 & 1 & 1.021 \\
    Osaka & 8,805,081 & 19 & 1.028 & 1.092 & 1.099 & 1.106 & 1.190 \\
    Hyogo & 5,550,574 & 12 & 1.026 & 1.075 & 1.081 & 1.088 & 1.154 \\
    Nara & 1,442,795 & 4 & 1 & 1 & 1 & 1 & 0.975 \\
    Wakayama & 1,069,912 & 3 & 1 & 1 & 1 & 1 & 0.964 \\
    Tottori & 613,289 & 2 & 1 & 1 & 1 & 1 & 0.829 \\
    Shimane & 761,503 & 2 & 1 & 1 & 1 & 1 & 1 \\
    Okayama & 1,950,828 & 5 & 1 & 1 & 1 & 1 & 1 \\
    Hiroshima & 2,878,915 & 7 & 1 & 1 & 1 & 1 & 1 \\
    Yamaguchi & 1,527,964 & 4 & 1 & 1 & 1 & 1 & 1 \\
    Tokushima & 824,108 & 3 & 0.914 & 0.830 & 0.819 & 0.807 & 0.742 \\
    Kagawa & 1,022,890 & 3 & 1 & 1 & 1 & 1 & 0.921 \\
    Ehime & 1,493,092 & 4 & 1 & 1 & 1 & 1 & 1 \\
    Kochi & 813,949 & 3 & 0.902 & 0.819 & 0.808 & 0.797 & 0.733 \\
    Fukuoka & 5,015,699 & 11 & 1.011 & 1.056 & 1.062 & 1.068 & 1.130 \\
    Saga & 876,654 & 3 & 0.972 & 0.883 & 0.871 & 0.859 & 0.790 \\
    Nagasaki & 1,516,523 & 4 & 1 & 1 & 1 & 1 & 1 \\
    Kumamoto & 1,859,344 & 5 & 1 & 1 & 1 & 1 & 1 \\
    Oita & 1,221,140 & 3 & 1 & 1 & 1 & 1 & 1 \\
    Miyazaki & 1,170,007 & 3 & 1 & 1 & 1 & 1 & 1 \\
    Kagoshima & 1,786,194 & 5 & 0.990 & 0.971 & 0.972 & 0.972 & 0.965 \\
    Okinawa & 1,318,220 & 4 & 0.974 & 0.929 & 0.926 & 0.923 & 0.891 \\
    \bottomrule
  \end{tabular}
\end{table}

\begin{table}[htbp]
  \centering
  \footnotesize
  \caption{Optimal weights $w_i^*$ of the prefectures for the 2009 election
($K=300$).}
  \label{tab:w2009}
  \begin{tabular}{lrrrrrrr}
    \toprule
    \textbf{Prefecture} & \textbf{$p_i$} & \textbf{$\bar{x}_i$} & \textbf{Adams} & \textbf{Dean} & \textbf{Hill} & \textbf{Webster} & \textbf{Jefferson} \\
    \midrule
    Hokkaido & 5,627,737 & 12 & 1.032 & 1.082 & 1.089 & 1.097 & 1.168 \\
    Aomori & 1,436,657 & 4 & 1 & 1 & 1 & 1 & 0.969 \\
    Iwate & 1,385,041 & 4 & 1 & 0.970 & 0.966 & 0.964 & 0.934 \\
    Miyagi & 2,360,218 & 6 & 1 & 1 & 1 & 1 & 1 \\
    Akita & 1,145,501 & 3 & 1 & 1 & 1 & 1 & 1 \\
    Yamagata & 1,216,181 & 3 & 1 & 1 & 1 & 1 & 1 \\
    Fukushima & 2,091,319 & 5 & 1 & 1 & 1 & 1 & 1 \\
    Ibaraki & 2,975,167 & 7 & 1 & 1 & 1 & 1 & 1.003 \\
    Tochigi & 2,016,631 & 5 & 1 & 1 & 1 & 1 & 1 \\
    Gunma & 2,024,135 & 5 & 1 & 1 & 1 & 1 & 1 \\
    Saitama & 7,054,243 & 15 & 1.035 & 1.093 & 1.100 & 1.109 & 1.189 \\
    Chiba & 6,056,462 & 13 & 1.025 & 1.078 & 1.085 & 1.093 & 1.167 \\
    Tokyo & 12,576,601 & 25 & 1.107 & 1.184 & 1.192 & 1.202 & 1.305 \\
    Kanagawa & 8,791,597 & 18 & 1.075 & 1.141 & 1.149 & 1.158 & 1.248 \\
    Niigata & 2,431,459 & 6 & 1 & 1 & 1 & 1 & 1 \\
    Toyama & 1,111,729 & 3 & 1 & 1 & 1 & 1 & 1.000 \\
    Ishikawa & 1,174,026 & 3 & 1 & 1 & 1 & 1 & 1 \\
    Fukui & 821,592 & 3 & 0.904 & 0.822 & 0.811 & 0.801 & 0.739 \\
    Yamanashi & 884,515 & 3 & 0.973 & 0.885 & 0.873 & 0.862 & 0.795 \\
    Nagano & 2,196,114 & 5 & 1 & 1 & 1 & 1 & 1 \\
    Gifu & 2,107,226 & 5 & 1 & 1 & 1 & 1 & 1 \\
    Shizuoka & 3,792,377 & 8 & 1.043 & 1.075 & 1.080 & 1.087 & 1.137 \\
    Aichi & 7,254,704 & 15 & 1.064 & 1.125 & 1.132 & 1.140 & 1.223 \\
    Mie & 1,866,963 & 5 & 1 & 1 & 1 & 1 & 1 \\
    Shiga & 1,380,361 & 4 & 1 & 0.966 & 0.963 & 0.961 & 0.931 \\
    Kyoto & 2,647,660 & 6 & 1 & 1 & 1 & 1 & 1.020 \\
    Osaka & 8,817,166 & 19 & 1.021 & 1.086 & 1.093 & 1.102 & 1.189 \\
    Hyogo & 5,590,601 & 12 & 1.025 & 1.075 & 1.082 & 1.090 & 1.160 \\
    Nara & 1,421,310 & 4 & 1 & 0.995 & 0.992 & 0.989 & 0.959 \\
    Wakayama & 1,035,969 & 3 & 1 & 1 & 1 & 1 & 0.932 \\
    Tottori & 607,012 & 2 & 1 & 1 & 1 & 0.986 & 0.819 \\
    Shimane & 742,223 & 2 & 1 & 1 & 1 & 1 & 1 \\
    Okayama & 1,957,264 & 5 & 1 & 1 & 1 & 1 & 1 \\
    Hiroshima & 2,876,642 & 7 & 1 & 1 & 1 & 1 & 1 \\
    Yamaguchi & 1,492,606 & 4 & 1 & 1 & 1 & 1 & 1 \\
    Tokushima & 809,950 & 3 & 0.891 & 0.810 & 0.799 & 0.789 & 0.728 \\
    Kagawa & 1,012,400 & 3 & 1 & 1 & 0.999 & 0.987 & 0.910 \\
    Ehime & 1,467,815 & 4 & 1 & 1 & 1 & 1 & 0.990 \\
    Kochi & 796,292 & 3 & 0.876 & 0.796 & 0.786 & 0.776 & 0.716 \\
    Fukuoka & 5,049,908 & 11 & 1.010 & 1.056 & 1.062 & 1.070 & 1.135 \\
    Saga & 866,369 & 3 & 0.953 & 0.866 & 0.855 & 0.844 & 0.779 \\
    Nagasaki & 1,478,632 & 4 & 1 & 1 & 1 & 1 & 0.997 \\
    Kumamoto & 1,842,233 & 5 & 1 & 0.995 & 0.996 & 0.997 & 0.994 \\
    Oita & 1,209,571 & 3 & 1 & 1 & 1 & 1 & 1 \\
    Miyazaki & 1,153,042 & 3 & 1 & 1 & 1 & 1 & 1 \\
    Kagoshima & 1,753,179 & 5 & 0.965 & 0.947 & 0.947 & 0.949 & 0.946 \\
    Okinawa & 1,361,594 & 4 & 0.999 & 0.953 & 0.950 & 0.948 & 0.918 \\
    \bottomrule
  \end{tabular}
\end{table}

\begin{table}[htbp]
  \centering
  \footnotesize
  \caption{Optimal weights $w_i^*$ of the prefectures for the 2012 election
($K=300$).}
  \label{tab:w2012}
  \begin{tabular}{lrrrrrrr}
    \toprule
    \textbf{Prefecture} & \textbf{$p_i$} & \textbf{$\bar{x}_i$} & \textbf{Adams} & \textbf{Dean} & \textbf{Hill} & \textbf{Webster} & \textbf{Jefferson} \\
    \midrule
    Hokkaido & 5,506,419 & 12 & 1.012 & 1.062 & 1.069 & 1.079 & 1.152 \\
    Aomori & 1,373,339 & 4 & 1 & 0.964 & 0.962 & 0.961 & 0.934 \\
    Iwate & 1,330,147 & 4 & 0.978 & 0.934 & 0.931 & 0.931 & 0.904 \\
    Miyagi & 2,348,165 & 6 & 1 & 1 & 1 & 1 & 1 \\
    Akita & 1,085,997 & 3 & 1 & 1 & 1 & 1 & 0.984 \\
    Yamagata & 1,168,924 & 3 & 1 & 1 & 1 & 1 & 1 \\
    Fukushima & 2,029,064 & 5 & 1 & 1 & 1 & 1 & 1 \\
    Ibaraki & 2,969,770 & 7 & 1 & 1 & 1 & 1 & 1.009 \\
    Tochigi & 2,007,683 & 5 & 1 & 1 & 1 & 1 & 1 \\
    Gunma & 2,008,068 & 5 & 1 & 1 & 1 & 1 & 1 \\
    Saitama & 7,194,556 & 15 & 1.058 & 1.118 & 1.126 & 1.137 & 1.223 \\
    Chiba & 6,216,289 & 13 & 1.055 & 1.110 & 1.118 & 1.128 & 1.207 \\
    Tokyo & 13,159,388 & 25 & 1.161 & 1.242 & 1.252 & 1.264 & 1.376 \\
    Kanagawa & 9,048,331 & 18 & 1.109 & 1.178 & 1.187 & 1.198 & 1.295 \\
    Niigata & 2,374,450 & 6 & 1 & 1 & 1 & 1 & 1 \\
    Toyama & 1,093,247 & 3 & 1 & 1 & 1 & 1 & 0.991 \\
    Ishikawa & 1,169,788 & 3 & 1 & 1 & 1 & 1 & 1 \\
    Fukui$^-$ & 806,314 & 3 & 0.889 & 0.808 & 0.798 & 0.790 & 0.731 \\
    Yamanashi$^-$ & 863,075 & 3 & 0.952 & 0.865 & 0.855 & 0.845 & 0.782 \\
    Nagano & 2,152,449 & 5 & 1 & 1 & 1 & 1 & 1 \\
    Gifu & 2,080,773 & 5 & 1 & 1 & 1 & 1 & 1 \\
    Shizuoka & 3,765,007 & 8 & 1.038 & 1.070 & 1.076 & 1.085 & 1.138 \\
    Aichi & 7,410,719 & 15 & 1.090 & 1.152 & 1.160 & 1.171 & 1.259 \\
    Mie & 1,854,724 & 5 & 1 & 1 & 1 & 1 & 1 \\
    Shiga & 1,410,777 & 4 & 1 & 0.990 & 0.988 & 0.987 & 0.959 \\
    Kyoto & 2,636,092 & 6 & 1 & 1 & 1 & 1 & 1.024 \\
    Osaka & 8,865,245 & 19 & 1.029 & 1.095 & 1.103 & 1.113 & 1.205 \\
    Hyogo & 5,588,133 & 12 & 1.027 & 1.077 & 1.085 & 1.095 & 1.169 \\
    Nara & 1,400,728 & 4 & 1 & 0.983 & 0.981 & 0.980 & 0.952 \\
    Wakayama & 1,002,198 & 3 & 1 & 1 & 0.992 & 0.982 & 0.908 \\
    Tottori & 588,667 & 2 & 1 & 1 & 1 & 0.961 & 0.800 \\
    Shimane & 717,397 & 2 & 1 & 1 & 1 & 1 & 0.975 \\
    Okayama & 1,945,276 & 5 & 1 & 1 & 1 & 1 & 1 \\
    Hiroshima & 2,860,750 & 7 & 1 & 1 & 1 & 1 & 1 \\
    Yamaguchi & 1,451,338 & 4 & 1 & 1 & 1 & 1 & 0.987 \\
    Tokushima$^-$ & 785,491 & 3 & 0.866 & 0.788 & 0.778 & 0.769 & 0.712 \\
    Kagawa & 995,842 & 3 & 1 & 0.998 & 0.986 & 0.975 & 0.903 \\
    Ehime & 1,431,493 & 4 & 1 & 1 & 1 & 1 & 0.973 \\
    Kochi$^-$ & 764,456 & 3 & 0.843 & 0.766 & 0.757 & 0.749 & 0.693 \\
    Fukuoka & 5,071,968 & 11 & 1.017 & 1.063 & 1.071 & 1.080 & 1.149 \\
    Saga$^-$ & 849,788 & 3 & 0.937 & 0.852 & 0.841 & 0.832 & 0.770 \\
    Nagasaki & 1,426,779 & 4 & 1 & 1 & 0.999 & 0.998 & 0.970 \\
    Kumamoto & 1,817,426 & 5 & 1 & 0.984 & 0.986 & 0.989 & 0.988 \\
    Oita & 1,196,529 & 3 & 1 & 1 & 1 & 1 & 1 \\
    Miyazaki & 1,135,233 & 3 & 1 & 1 & 1 & 1 & 1 \\
    Kagoshima & 1,706,242 & 5 & 0.941 & 0.924 & 0.925 & 0.928 & 0.928 \\
    Okinawa & 1,392,818 & 4 & 1 & 0.978 & 0.975 & 0.974 & 0.947 \\
    \bottomrule
  \end{tabular}
\end{table}

\begin{table}[htbp]
  \centering
  \footnotesize
  \caption{Optimal weights $w_i^*$ of the prefectures for the 2014 election
($K=295$).}
  \label{tab:w2014}
  \begin{tabular}{lrrrrrrr}
    \toprule
    \textbf{Prefecture} & \textbf{$p_i$} & \textbf{$\bar{x}_i$} & \textbf{Adams} & \textbf{Dean} & \textbf{Hill} & \textbf{Webster} & \textbf{Jefferson} \\
    \midrule
    Hokkaido & 5,506,419 & 12 & 1 & 1.020 & 1.026 & 1.034 & 1.100 \\
    Aomori$^-$ & 1,373,339 & 4 & 0.971 & 0.926 & 0.923 & 0.921 & 0.892 \\
    Iwate$^-$ & 1,330,147 & 4 & 0.941 & 0.897 & 0.894 & 0.892 & 0.864 \\
    Miyagi & 2,348,165 & 6 & 0.997 & 0.995 & 0.998 & 1 & 1 \\
    Akita & 1,085,997 & 3 & 1 & 1 & 1 & 1 & 0.940 \\
    Yamagata & 1,168,924 & 3 & 1 & 1 & 1 & 1 & 1 \\
    Fukushima & 2,029,064 & 5 & 1 & 1 & 1 & 1 & 1 \\
    Ibaraki & 2,969,770 & 7 & 1 & 1 & 1 & 1 & 1 \\
    Tochigi & 2,007,683 & 5 & 1 & 1 & 1 & 1 & 1 \\
    Gunma & 2,008,068 & 5 & 1 & 1 & 1 & 1 & 1 \\
    Saitama & 7,194,556 & 15 & 1.018 & 1.074 & 1.081 & 1.089 & 1.168 \\
    Chiba & 6,216,289 & 13 & 1.015 & 1.066 & 1.073 & 1.081 & 1.154 \\
    Tokyo & 13,159,388 & 25 & 1.117 & 1.193 & 1.202 & 1.211 & 1.315 \\
    Kanagawa & 9,048,331 & 18 & 1.067 & 1.131 & 1.139 & 1.148 & 1.237 \\
    Niigata & 2,374,450 & 6 & 1 & 1 & 1 & 1 & 1 \\
    Toyama & 1,093,247 & 3 & 1 & 1 & 1 & 1 & 0.947 \\
    Ishikawa & 1,169,788 & 3 & 1 & 1 & 1 & 1 & 1 \\
    Fukui & 806,314 & 2 & 1 & 1 & 1 & 1 & 1 \\
    Yamanashi & 863,075 & 2 & 1 & 1 & 1 & 1 & 1 \\
    Nagano & 2,152,449 & 5 & 1 & 1 & 1 & 1 & 1 \\
    Gifu & 2,080,773 & 5 & 1 & 1 & 1 & 1 & 1 \\
    Shizuoka & 3,765,007 & 8 & 1 & 1.027 & 1.033 & 1.040 & 1.087 \\
    Aichi & 7,410,719 & 15 & 1.048 & 1.106 & 1.114 & 1.122 & 1.203 \\
    Mie$^-$ & 1,854,724 & 5 & 0.984 & 0.964 & 0.966 & 0.967 & 0.964 \\
    Shiga & 1,410,777 & 4 & 0.998 & 0.951 & 0.948 & 0.946 & 0.916 \\
    Kyoto & 2,636,092 & 6 & 1 & 1 & 1 & 1 & 1 \\
    Osaka & 8,865,245 & 19 & 1 & 1.051 & 1.059 & 1.067 & 1.152 \\
    Hyogo & 5,588,133 & 12 & 1 & 1.035 & 1.042 & 1.049 & 1.117 \\
    Nara$^-$ & 1,400,728 & 4 & 0.991 & 0.944 & 0.941 & 0.939 & 0.910 \\
    Wakayama & 1,002,198 & 3 & 1 & 0.965 & 0.953 & 0.941 & 0.868 \\
    Tottori & 588,667 & 2 & 1 & 1 & 0.969 & 0.921 & 0.765 \\
    Shimane & 717,397 & 2 & 1 & 1 & 1 & 1 & 0.932 \\
    Okayama & 1,945,276 & 5 & 1 & 1 & 1 & 1 & 1 \\
    Hiroshima & 2,860,750 & 7 & 1 & 1 & 1 & 1 & 1 \\
    Yamaguchi & 1,451,338 & 4 & 1 & 0.978 & 0.975 & 0.973 & 0.943 \\
    Tokushima & 785,491 & 2 & 1 & 1 & 1 & 1 & 1 \\
    Kagawa & 995,842 & 3 & 1 & 0.959 & 0.946 & 0.935 & 0.862 \\
    Ehime & 1,431,493 & 4 & 1 & 0.965 & 0.962 & 0.960 & 0.930 \\
    Kochi & 764,456 & 2 & 1 & 1 & 1 & 1 & 0.993 \\
    Fukuoka & 5,071,968 & 11 & 1 & 1.021 & 1.028 & 1.035 & 1.098 \\
    Saga & 849,788 & 2 & 1 & 1 & 1 & 1 & 1 \\
    Nagasaki & 1,426,779 & 4 & 1 & 0.962 & 0.959 & 0.957 & 0.927 \\
    Kumamoto$^-$ & 1,817,426 & 5 & 0.964 & 0.945 & 0.946 & 0.948 & 0.944 \\
    Oita & 1,196,529 & 3 & 1 & 1 & 1 & 1 & 1 \\
    Miyazaki & 1,135,233 & 3 & 1 & 1 & 1 & 1 & 0.983 \\
    Kagoshima$^-$ & 1,706,242 & 5 & 0.905 & 0.887 & 0.888 & 0.890 & 0.887 \\
    Okinawa & 1,392,818 & 4 & 0.985 & 0.939 & 0.936 & 0.934 & 0.905 \\
    \bottomrule
  \end{tabular}
\end{table}

\begin{table}[htbp]
  \centering
  \footnotesize
  \caption{Optimal weights $w_i^*$ of the prefectures for the 2017 election
($K=289$).}
  \label{tab:w2017}
  \begin{tabular}{lrrrrrrr}
    \toprule
    \textbf{Prefecture} & \textbf{$p_i$} & \textbf{$\bar{x}_i$} & \textbf{Adams} & \textbf{Dean} & \textbf{Hill} & \textbf{Webster} & \textbf{Jefferson} \\
    \midrule
    Hokkaido & 5,348,768 & 12 & 1 & 1 & 1.005 & 1.012 & 1.080 \\
    Aomori & 1,302,132 & 3 & 1 & 1 & 1 & 1 & 1 \\
    Iwate & 1,272,745 & 3 & 1 & 1 & 1 & 1 & 1 \\
    Miyagi & 2,291,508 & 6 & 0.972 & 0.977 & 0.981 & 0.986 & 1 \\
    Akita & 1,017,149 & 3 & 1 & 0.986 & 0.974 & 0.962 & 0.890 \\
    Yamagata & 1,116,752 & 3 & 1 & 1 & 1 & 1 & 0.977 \\
    Fukushima & 1,898,880 & 5 & 1 & 0.994 & 0.996 & 0.998 & 0.997 \\
    Ibaraki & 2,862,997 & 7 & 1 & 1 & 1 & 1 & 1 \\
    Tochigi & 1,927,885 & 5 & 1 & 1 & 1 & 1 & 1 \\
    Gunma & 1,930,380 & 5 & 1 & 1 & 1 & 1 & 1 \\
    Saitama & 7,111,168 & 15 & 1.005 & 1.068 & 1.077 & 1.085 & 1.167 \\
    Chiba & 6,047,216 & 13 & 1 & 1.043 & 1.052 & 1.060 & 1.134 \\
    Tokyo & 12,948,463 & 25 & 1.098 & 1.181 & 1.191 & 1.201 & 1.308 \\
    Kanagawa & 8,887,304 & 18 & 1.047 & 1.118 & 1.127 & 1.136 & 1.228 \\
    Niigata & 2,289,345 & 6 & 0.971 & 0.976 & 0.981 & 0.985 & 1 \\
    Toyama & 1,052,353 & 3 & 1 & 1 & 1 & 0.996 & 0.921 \\
    Ishikawa & 1,138,322 & 3 & 1 & 1 & 1 & 1 & 0.996 \\
    Fukui & 774,337 & 2 & 1 & 1 & 1 & 1 & 1 \\
    Yamanashi & 819,205 & 2 & 1 & 1 & 1 & 1 & 1 \\
    Nagano & 2,067,713 & 5 & 1 & 1 & 1 & 1 & 1 \\
    Gifu & 1,989,980 & 5 & 1 & 1 & 1 & 1 & 1 \\
    Shizuoka & 3,626,991 & 8 & 1 & 1 & 1.003 & 1.009 & 1.058 \\
    Aichi & 7,260,847 & 15 & 1.027 & 1.091 & 1.100 & 1.108 & 1.192 \\
    Mie & 1,776,805 & 4 & 1 & 1 & 1 & 1 & 1 \\
    Shiga & 1,386,795 & 4 & 0.980 & 0.941 & 0.939 & 0.937 & 0.910 \\
    Kyoto & 2,533,645 & 6 & 1 & 1 & 1 & 1 & 1 \\
    Osaka & 8,524,530 & 19 & 1 & 1.017 & 1.026 & 1.034 & 1.119 \\
    Hyogo & 5,398,880 & 12 & 1 & 1.006 & 1.014 & 1.022 & 1.090 \\
    Nara & 1,351,535 & 3 & 1 & 1 & 1 & 1 & 1 \\
    Wakayama & 956,199 & 3 & 1 & 0.927 & 0.916 & 0.905 & 0.837 \\
    Tottori & 567,993 & 2 & 1 & 0.991 & 0.942 & 0.896 & 0.746 \\
    Shimane & 687,180 & 2 & 1 & 1 & 1 & 1 & 0.902 \\
    Okayama & 1,885,691 & 5 & 1.000 & 0.987 & 0.989 & 0.991 & 0.990 \\
    Hiroshima & 2,795,626 & 7 & 0.988 & 1 & 1 & 1 & 1 \\
    Yamaguchi & 1,390,689 & 4 & 0.983 & 0.943 & 0.942 & 0.940 & 0.913 \\
    Tokushima & 747,141 & 2 & 1 & 1 & 1 & 1 & 0.981 \\
    Kagawa & 961,844 & 3 & 1 & 0.932 & 0.921 & 0.910 & 0.842 \\
    Ehime & 1,365,508 & 4 & 0.965 & 0.926 & 0.925 & 0.923 & 0.896 \\
    Kochi & 722,728 & 2 & 1 & 1 & 1 & 1 & 0.949 \\
    Fukuoka & 4,995,297 & 11 & 1 & 1.012 & 1.020 & 1.028 & 1.093 \\
    Saga & 827,702 & 2 & 1 & 1 & 1 & 1 & 1 \\
    Nagasaki & 1,365,241 & 4 & 0.965 & 0.926 & 0.925 & 0.923 & 0.896 \\
    Kumamoto & 1,771,440 & 4 & 1 & 1 & 1 & 1 & 1 \\
    Oita & 1,150,436 & 3 & 1 & 1 & 1 & 1 & 1 \\
    Miyazaki & 1,096,407 & 3 & 1 & 1 & 1 & 1 & 0.960 \\
    Kagoshima & 1,631,662 & 4 & 1 & 1 & 1 & 1 & 1 \\
    Okinawa & 1,410,487 & 4 & 0.997 & 0.957 & 0.955 & 0.953 & 0.926 \\
    \bottomrule
  \end{tabular}
\end{table}

\begin{table}[htbp]
  \centering
  \footnotesize
  \caption{Optimal weights $w_i^*$ of the prefectures for the 2021 election
($K=289$).}
  \label{tab:w2021}
  \begin{tabular}{lrrrrrrr}
    \toprule
    \textbf{Prefecture} & \textbf{$p_i$} & \textbf{$\bar{x}_i$} & \textbf{Adams} & \textbf{Dean} & \textbf{Hill} & \textbf{Webster} & \textbf{Jefferson} \\
    \midrule
    Hokkaido & 5,151,366 & 12 & 1 & 1 & 1 & 1.003 & 1.075 \\
    Aomori & 1,224,334 & 3 & 1 & 1 & 1 & 1 & 1 \\
    Iwate & 1,194,745 & 3 & 1 & 1 & 1 & 1 & 1 \\
    Miyagi$^-$ & 2,242,701 & 6 & 0.978 & 0.983 & 0.988 & 0.992 & 1 \\
    Akita & 950,192 & 3 & 1 & 0.947 & 0.936 & 0.925 & 0.859 \\
    Yamagata & 1,056,617 & 3 & 1 & 1 & 1 & 1 & 0.956 \\
    Fukushima$^-$ & 1,797,450 & 5 & 0.979 & 0.967 & 0.970 & 0.972 & 0.975 \\
    Ibaraki & 2,763,432 & 7 & 1 & 1 & 1 & 1 & 1 \\
    Tochigi & 1,869,256 & 5 & 1 & 1 & 1 & 1 & 1 \\
    Gunma & 1,860,213 & 5 & 1 & 1 & 1 & 1 & 1 \\
    Saitama$^+$ & 7,055,679 & 15 & 1.025 & 1.089 & 1.099 & 1.108 & 1.196 \\
    Chiba$^+$ & 6,040,814 & 13 & 1.013 & 1.071 & 1.080 & 1.089 & 1.171 \\
    Tokyo$^+$ & 13,233,213 & 25 & 1.154 & 1.241 & 1.252 & 1.263 & 1.381 \\
    Kanagawa$^+$ & 8,876,834 & 18 & 1.075 & 1.148 & 1.158 & 1.168 & 1.268 \\
    Niigata$^-$ & 2,168,416 & 6 & 0.945 & 0.950 & 0.955 & 0.959 & 0.981 \\
    Toyama & 1,009,388 & 3 & 1 & 1 & 0.994 & 0.983 & 0.913 \\
    Ishikawa & 1,102,692 & 3 & 1 & 1 & 1 & 1 & 0.997 \\
    Fukui & 746,068 & 2 & 1 & 1 & 1 & 1 & 1 \\
    Yamanashi & 786,673 & 2 & 1 & 1 & 1 & 1 & 1 \\
    Nagano & 1,997,083 & 5 & 1 & 1 & 1 & 1 & 1 \\
    Gifu & 1,900,727 & 5 & 1 & 1 & 1 & 1 & 1 \\
    Shizuoka & 3,523,817 & 8 & 1 & 1 & 1.002 & 1.009 & 1.062 \\
    Aichi$^+$ & 7,164,758 & 15 & 1.041 & 1.106 & 1.116 & 1.125 & 1.215 \\
    Mie & 1,698,876 & 4 & 1 & 1 & 1 & 1 & 1 \\
    Shiga$^-$ & 1,359,501 & 4 & 0.988 & 0.948 & 0.947 & 0.945 & 0.922 \\
    Kyoto & 2,460,764 & 6 & 1 & 1 & 1 & 1 & 1 \\
    Osaka & 8,411,856 & 19 & 1 & 1.032 & 1.041 & 1.050 & 1.141 \\
    Hyogo & 5,215,927 & 12 & 1 & 1 & 1.007 & 1.015 & 1.089 \\
    Nara & 1,302,289 & 3 & 1 & 1 & 1 & 1 & 1 \\
    Wakayama$^-$ & 911,408 & 3 & 0.993 & 0.908 & 0.898 & 0.887 & 0.824 \\
    Tottori & 543,979 & 2 & 1 & 0.975 & 0.928 & 0.883 & 0.738 \\
    Shimane & 656,197 & 2 & 1 & 1 & 1 & 1 & 0.890 \\
    Okayama$^-$ & 1,810,900 & 5 & 0.987 & 0.974 & 0.977 & 0.979 & 0.983 \\
    Hiroshima$^-$ & 2,713,163 & 7 & 0.986 & 1 & 1 & 1 & 1 \\
    Yamaguchi$^-$ & 1,313,593 & 4 & 0.954 & 0.916 & 0.915 & 0.913 & 0.891 \\
    Tokushima & 693,084 & 2 & 1 & 1 & 1 & 1 & 0.940 \\
    Kagawa & 920,884 & 3 & 1 & 0.917 & 0.907 & 0.896 & 0.833 \\
    Ehime$^-$ & 1,296,806 & 4 & 0.942 & 0.904 & 0.903 & 0.902 & 0.880 \\
    Kochi & 677,762 & 2 & 1 & 1 & 1 & 1 & 0.919 \\
    Fukuoka & 4,923,076 & 11 & 1 & 1.025 & 1.034 & 1.042 & 1.113 \\
    Saga & 796,285 & 2 & 1 & 1 & 1 & 1 & 1 \\
    Nagasaki$^-$ & 1,296,014 & 4 & 0.942 & 0.904 & 0.903 & 0.901 & 0.879 \\
    Kumamoto & 1,705,205 & 4 & 1 & 1 & 1 & 1 & 1 \\
    Oita & 1,098,735 & 3 & 1 & 1 & 1 & 1 & 0.994 \\
    Miyazaki & 1,047,871 & 3 & 1 & 1 & 1 & 1 & 0.948 \\
    Kagoshima & 1,545,763 & 4 & 1 & 1 & 1 & 1 & 1 \\
    Okinawa & 1,424,749 & 4 & 1 & 0.993 & 0.992 & 0.991 & 0.966 \\
    \bottomrule
  \end{tabular}
\end{table}

\begin{table}[htbp]
  \centering
  \footnotesize
  \caption{Optimal weights $w_i^*$ of the prefectures for the 2024 election
($K=289$).}
  \label{tab:w2024}
  \begin{tabular}{lrrrrrrr}
    \toprule
    \textbf{Prefecture} & \textbf{$p_i$} & \textbf{$\bar{x}_i$} & \textbf{Adams} & \textbf{Dean} & \textbf{Hill} & \textbf{Webster} & \textbf{Jefferson} \\
    \midrule
    Hokkaido & 5,151,366 & 12 & 1 & 1 & 1 & 1.001 & 1.068 \\
    Aomori & 1,224,334 & 3 & 1 & 1 & 1 & 1 & 1 \\
    Iwate & 1,194,745 & 3 & 1 & 1 & 1 & 1 & 1 \\
    Miyagi & 2,242,701 & 5 & 1 & 1 & 1 & 1 & 1.007 \\
    Akita & 950,192 & 3 & 1 & 0.944 & 0.934 & 0.923 & 0.853 \\
    Yamagata & 1,056,617 & 3 & 1 & 1 & 1 & 1 & 0.949 \\
    Fukushima & 1,797,450 & 4 & 1 & 1 & 1 & 1 & 1 \\
    Ibaraki & 2,763,432 & 7 & 1 & 1 & 1 & 1 & 1 \\
    Tochigi & 1,869,256 & 5 & 1 & 1 & 1 & 1 & 1 \\
    Gunma & 1,860,213 & 5 & 1 & 0.998 & 1 & 1 & 1 \\
    Saitama & 7,055,679 & 16 & 1 & 1.021 & 1.030 & 1.039 & 1.118 \\
    Chiba & 6,040,814 & 14 & 1 & 1 & 1.004 & 1.012 & 1.085 \\
    Tokyo & 13,233,213 & 30 & 1 & 1.035 & 1.045 & 1.054 & 1.150 \\
    Kanagawa & 8,876,834 & 20 & 1 & 1.033 & 1.043 & 1.052 & 1.139 \\
    Niigata & 2,168,416 & 5 & 1 & 1 & 1 & 1 & 1 \\
    Toyama & 1,009,388 & 3 & 1 & 1 & 0.992 & 0.981 & 0.907 \\
    Ishikawa & 1,102,692 & 3 & 1 & 1 & 1 & 1 & 0.990 \\
    Fukui & 746,068 & 2 & 1 & 1 & 1 & 1 & 1 \\
    Yamanashi & 786,673 & 2 & 1 & 1 & 1 & 1 & 1 \\
    Nagano & 1,997,083 & 5 & 1 & 1 & 1 & 1 & 1 \\
    Gifu & 1,900,727 & 5 & 1 & 1 & 1 & 1 & 1 \\
    Shizuoka & 3,523,817 & 8 & 1 & 1 & 1 & 1.007 & 1.055 \\
    Aichi & 7,164,758 & 16 & 1 & 1.036 & 1.046 & 1.055 & 1.136 \\
    Mie & 1,698,876 & 4 & 1 & 1 & 1 & 1 & 1 \\
    Shiga & 1,359,501 & 3 & 1 & 1 & 1 & 1 & 1 \\
    Kyoto & 2,460,764 & 6 & 1 & 1 & 1 & 1 & 1 \\
    Osaka & 8,411,856 & 19 & 1 & 1.029 & 1.039 & 1.048 & 1.133 \\
    Hyogo & 5,215,927 & 12 & 1 & 1 & 1.006 & 1.013 & 1.081 \\
    Nara & 1,302,289 & 3 & 1 & 1 & 1 & 1 & 1 \\
    Wakayama & 911,408 & 2 & 1 & 1 & 1 & 1 & 1 \\
    Tottori & 543,979 & 2 & 1 & 0.973 & 0.926 & 0.881 & 0.733 \\
    Shimane & 656,197 & 2 & 1 & 1 & 1 & 1 & 0.884 \\
    Okayama & 1,810,900 & 4 & 1 & 1 & 1 & 1 & 1 \\
    Hiroshima & 2,713,163 & 6 & 1 & 1.001 & 1.008 & 1.014 & 1.044 \\
    Yamaguchi & 1,313,593 & 3 & 1 & 1 & 1 & 1 & 1 \\
    Tokushima & 693,084 & 2 & 1 & 1 & 1 & 1 & 0.934 \\
    Kagawa & 920,884 & 3 & 1 & 0.915 & 0.905 & 0.895 & 0.827 \\
    Ehime & 1,296,806 & 3 & 1 & 1 & 1 & 1 & 1 \\
    Kochi & 677,762 & 2 & 1 & 1 & 1 & 1 & 0.913 \\
    Fukuoka & 4,923,076 & 11 & 1 & 1.023 & 1.032 & 1.040 & 1.105 \\
    Saga & 796,285 & 2 & 1 & 1 & 1 & 1 & 1 \\
    Nagasaki & 1,296,014 & 3 & 1 & 1 & 1 & 1 & 1 \\
    Kumamoto & 1,705,205 & 4 & 1 & 1 & 1 & 1 & 1 \\
    Oita & 1,098,735 & 3 & 1 & 1 & 1 & 1 & 0.987 \\
    Miyazaki & 1,047,871 & 3 & 1 & 1 & 1 & 1 & 0.941 \\
    Kagoshima & 1,545,763 & 4 & 1 & 1 & 1 & 1 & 1 \\
    Okinawa & 1,424,749 & 4 & 1 & 0.991 & 0.990 & 0.989 & 0.960 \\
    \bottomrule
  \end{tabular}
\end{table}

\end{document}